\documentclass[11pt]{article}

\usepackage{amsfonts,amssymb,amsmath,latexsym,indentfirst}

\newtheorem{prop}{Proposition}[section]
\newtheorem{rema}{Remark}[section]
\newtheorem{defi}{Definition}[section]
\newtheorem{lemm}{Lemma}[section]
\newtheorem{theo}{Theorem}[section]

\newcommand{\R}{\ensuremath{{\mathbb{R}} }}

\newcommand{\peq}{\hspace*{0.10in}}
\newcommand{\peqq}{\hspace*{0.05in}}

\newcommand{\quebra}{\hspace*{-0.25in}}

\newcommand{\fim}{\rightline{$\blacksquare$}}

\author{Luiz Gustavo  Farah and Felipe Linares\\
\\
IMPA\\
Estrada Dona Castorina 110,\\
Rio de Janeiro 22460-320, Brazil}
\title{Global rough solutions to the cubic nonlinear Boussinesq equation
\footnotetext{The authors were partially supported by CNPq-Brazil.}}
\date{}

\begin{document}
\maketitle

\begin{abstract}
We prove that the initial value problem (IVP) for the cubic defocusing nonlinear Boussinesq
equation $u_{tt}-u_{xx}+u_{xxxx}-(|u|^2u)_{xx}=0$ on the real line is globally well-posed in
$H^{s}(\R)$ provided $2/3<s<1$.
\end{abstract}

\section{Introduction}
In this work we consider the initial value problem (IVP) for the
Boussi\-nesq-type equation
\begin{eqnarray}\label{NLB}
\left\{
\begin{array}{l}
u_{tt}-u_{xx}+u_{xxxx}-(|u|^2u)_{xx}=0, \peq  x\in \R,\, t>0,\\
u(x,0)=\phi(x); \peqq \partial_tu(x,0)=\psi_x(x).
\end{array} \right.
\end{eqnarray}

Equations of this type are generalization of the one originally derived by Boussinesq \cite{BOU}
in his study of nonlinear, dispersive wave propagation. The equation \eqref{NLB} has also been
used as a model of nonlinear strings by Zakharov \cite{Z} and in the study of shape-memory alloys
by Falk \textit{et al} \cite{FLS}.

It is know that equation (\ref{NLB}) is locally well-posed in time in $H^1(\mathbb R)$ (see Linares \cite{FL}).
Moreover, this local solution has the following conserved energy,
\begin{equation}\label{EC}
E(u)(t)=\frac{1}{2}\|u(t)\|^2_{H^1}+ \frac{1}{2}\|(-\Delta)^{-\frac{1}{2}}\partial_tu(t)\|^2_{L^2}
+ \frac{1}{4}\|u(t)\|^4_{L^4}.
\end{equation}

Local theory together with this conserved quantity immediately yield global-in-time well-posedness
of (\ref{NLB}) from data $(\phi, \psi)\in H^1(\R)\times L^2(\R)$.

Our principal aim is to loosen the regularity requirements on the
initial data which ensure global-in-time solutions for the IVP
(\ref{NLB}). Here we use the approach introduced by Colliander,
Keel, Staffilani, Takaoka and Tao in \cite{CKSTT4}, called the
I-method. Since the equation (\ref{NLB}) has two derivatives in time, it is not clear whether the refined approach introduced by the same authors in \cite{CKSTT6} and \cite{CKSTT5} can be use to improve our global result stated in Theorem \ref{T1} below.

Note that when $(\phi, \psi)\in H^{s}(\R)\times H^{s-1}(\R)$ with
$s<1$ in (\ref{NLB}), the energy could be infinite, and so the
conservation law (\ref{EC}) is meaningless. To overcome this
difficulty, we follow the $I$-method scheme and introduce a modified
energy functional which is also defined for less regular functions.
Unfortunately, this new functional is not strictly conserved, but we
can show that it is \textit{almost} conserved in time. When one is
able to control its growth in time explicitly this allows to iterate
a modified local existence theorem to continue the solution to any
time $T$. This method has been successfully applied to
several equations which have a scaling invariance with sometimes
even optimal global well-posedness results. Although for the cubic
Boussinesq equation (\ref{NLB}) such a scaling argument does not work and there is no conservation law at level $L^2$, we can also apply
this method to prove global existence result with rough initial data.

We shall notice that the local well-posedness theory for the IVP \eqref{NLB} with initial data
$(\phi, \psi)\in H^{s}(\R)\times H^{s-1}(\R)$, $0<s<1$, seems to be, up to our knowledge, unknown
in the literature. (The case $s=0$ was obtained in \cite{FL}). Thus we will also establish it to
develop our analysis. This will be done using the ideas introduced by Farah \cite{LG4} in his study
of the \lq\lq good'' Boussinesq equation.

To describe our well-posedness results we define next the $X_{s,b}$ spaces related to our problem. These
spaces were first defined by Fang and Grillakis \cite{FG} for the Boussinesq-type
equations in the periodic setting. Using these spaces and following Bourgain's argument introduced in
\cite{B}
they proved local well-posedness for (\ref{NLB}) with the spatial variable in the unit circle assuming
$u_0 \in H^{s}(\mathbb T)$, $u_1 \in H^{-2+s}(\mathbb T)$, with $0\leq s\leq 1$ and  $|f(u)|\leq c|u|^{p}$,
with
$1<p<\frac{3-2s}{1-2s}$ if $0\leq s<\frac{1}{2}$ and $1<p<\infty$ if $\frac{1}{2}\leq s\leq 1$. Moreover,
if $u_0 \in H^{1}(\mathbb T)$, $u_1 \in H^{-1}(\mathbb T)$ and  $f(u)= \lambda|u|^{q-1}u-|u|^{p-1}u$, with
$1<q<p$ and $\lambda \in \R$ then the solution is global.

Next we give the precise definition of the $X_{s,b}$ spaces for the Boussinesq-type equation in the
continuous case.

\begin{defi}\label{BOS}
For $s,b \in \R$, $X_{s,b}$ denotes the completion of the Schwartz class $\mathcal{S}(\R^2)$ with respect
to the norm
\begin{equation*}
\|F\|_{X_{s,b}}=\|\langle|\tau|-\gamma(\xi)\rangle^b\langle\xi\rangle^s \widetilde{F}\|_{L^{2}_{\xi,\tau}}
\end{equation*}
where $\gamma(\xi)\equiv\sqrt{{\xi}^2+{\xi}^4}$, $\sim$ denotes the time-space Fourier transform and
$\langle a\rangle\equiv (1+|a|^2)^{\frac{1}{2}}$.
\end{defi}

We will also need the localized $X_{s,b}$ spaces defined as follows.

\begin{defi}\label{BL}
For $s,b \in \R$ and $\delta \geq 0$, $X_{s,b}^{\delta}$ denotes the space endowed with the norm
\begin{equation*}
\|u\|_{X_{s,b}^{\delta}}=\inf_{w\in X_{s,b}}\left\{\|w\|_{X_{s,b}}:w(t)=u(t) \textrm{ on }  [0,{\delta}]\right\}.
\end{equation*}
\end{defi}

Now we state the main results of this paper.

\begin{theo}\label{t3.1}
Let $s\geq 0$, then for all $\phi\in H^s(\R)$ and $\psi\in H^{s-1}(\R)$, there exist $\delta=\delta(\|\phi\|_{H^s},\|\psi\|_{H^{s-1}})$ and a unique solution $u$ of the IVP (\ref{NLB})
\begin{equation*}
u\in C([0,\delta]:H^s(\R))\cap X^{\delta}_{s,b}.
\end{equation*}

Moreover, given $\delta'\in (0,\delta)$ there exists $R=R(\delta')>0$ such that giving the set $ W\equiv\{(\tilde{\phi},\tilde{\psi})\in H^s(\R)\times H^{s-1}(\R):\|\tilde{\phi}-\phi\|_{H^s(\R)}^2+ \|\tilde{\psi}-\psi\|_{H^{s-1}(\R)}^2<R\}$ the map solution
\begin{equation*}
S:W \longrightarrow C([0,\delta']:H^s(\R))\cap X^{\delta}_{s,b}, \peq (\tilde{\phi},\tilde{\psi})\longmapsto u(t)
\end{equation*}
is Lipschitz.
\end{theo}

\begin{rema}\label{obs1}  The same result holds when we consider the equation \eqref{NLB} with positive sign
in front of the nonlinearity.
\end{rema}

\begin{theo}\label{T1}
The initial value problem (\ref{NLB}) is globally well-posed in $H^s(\R)$ for all \peqq $2/3<s<1$. Moreover the solution satisfies
\begin{equation}\label{pb}
\sup_{t\in[0,T]}\left\{\|u(t)\|^2_{H^{s}}+ \|(-\Delta)^{-\frac{1}{2}}\partial_tu(t)\|^2_{H^{s-1}}\right\}\leq C(1+T)^{\frac{1-s}{6s-4}+}
\end{equation}
where the constant $C$ depends only on $s$, $\|\phi\|_{H^{s}}$ and $\|\psi\|_{H^{s-1}}$.
\end{theo}

\begin{rema} Solutions of the IVP associated to the equation in Remark \ref{obs1} may blow-up
for arbitrary initial data, see for instance Liu \cite{LIU2} and Angulo and Scialom \cite{AS}.
\end{rema}

We observe that the $X_{s,b}$ spaces used here are very similar to the Bourgain
spaces associated to the Schr\"odinger equation (see Section \ref{S3} below). This
allows us to follow some ideas employed to study well posedness for the IVP associated to
the cubic Schr\"o\-din\-ger equation.
However the global problem is quite different in each case. To illustrate this we consider
the IVP associated to the cubic Schr\"o\-din\-ger equation, i.e.,
\begin{eqnarray}\label{NLS}
\left\{
\begin{array}{l}
i\phi_{t}+\phi_{xx}=|\phi|^2\phi, \peq  x\in \R^n,\, t>0,\\
\phi(x,0)=\phi_0(x).
\end{array} \right.
\end{eqnarray}

Solutions of the IVP \eqref{NLS} leave the $L^2$-norm invariant, i.e.,
\begin{equation}\label{MC}
\|\phi(\cdot,t)\|_{L^2(\R^n)}=\|\phi_0\|_{L^2(\R^n)}
\end{equation}
and the same holds for the energy, i.e.,
\begin{equation}\label{ESC}
E(\phi)(t)=\frac{1}{2}\|\nabla_x\phi(x,t)\|^2_{L^2(\R^n)}+ \frac{1}{4}\|\phi(x,t)\|^4_{L^4(\R^n)}.
\end{equation}

It is well known that the equation \eqref{NLS} is
locally well-posed in $L^2(\R)$ (resp. $H^1(\R)$) in the sub-critical sense, that is, the existence time
$T$ given in the local theory depends only on the norm of the initial data. Therefore, mass conservation
(\ref{MC}) (resp. energy conservation (\ref{ESC})) immediately yields global well-posedness in $L^2(\R)$
(resp. $H^1(\R)$). Hence, it is natural to expect that the IVP (\ref{NLS}) is globally well-posed in
$H^s(\R)$, with $0<s<1$.

The situation for the cubic Boussinesq equation (\ref{NLB}) is much more different. Mass conservation
does not hold for this equation. In this case, we cannot apply the same argument used for the sub-critical
Schr\"odinger equation (\ref{NLS}) to obtain global solutions for initial data
$(\phi, \psi)\in L^{2}(\R)\times H^{-1}(\R)$. In fact, it is an open problem to prove global well-posedness
in this case.

We should notice, that this is similar to the case analyzed by J. Colliander, M. Keel, G. Staffilani, H.
Takaoka and T. Tao \cite{CKSTT4}, where much of our calculations are based. In \cite{CKSTT4}, the authors
study, in particular, the cubic Schr\"odinger equation (\ref{NLS}), with the space variable belonging to
$\R^2$. Such equation is know to be critical and it is also an open problem to prove global well-posedness
in $L^2(\R^2)$.

The plan of this paper is as follows. In the next section we derive the modified energy functional.
Next we introduce some notation and prove refined Strichartz estimates for the solutions of the Boussinesq
equation. In Section 4, we study the local theory related to the IVP (\ref{NLB}). Finally, in Sections 5 and
6, we prove the almost conserved law, which implies the global result stated in Theorem \ref{T1}.

\section{Modified energy functional}

As we mention in the introduction, we follow the $I$-method scheme and introduce a smoothed version of
solutions of (\ref{NLB}) with finite energy and show that its  energy is \textit{almost} conserved
in time. This smoothed version is expressed as follows.

Given $s<1$ and a parameter $N\gg 1$, define the multiplier operator
\begin{equation*}
\widehat{I_Nf(\xi)}\equiv m_N(\xi)\widehat{f}(\xi),
\end{equation*}
where the multiplier $m_N(\xi)$ is smooth, radially symmetric, nondecreasing in $|\xi|$ and
\begin{eqnarray*}
m_N(\xi)=\left\{
\begin{array}{l l }
1&, \textrm{ if } |\xi|\leq N,\\
\left(\dfrac{N}{|\xi|}\right)^{1-s}&, \textrm{ if } |\xi|\geq N.
\end{array} \right.
\end{eqnarray*}

To simplify the notation, we omit the dependence of $N$ in $I_N$ and denote it only by $I$. Note that the operator $I$ is smooth of order $1-s$. Indeed, we have
\begin{equation}\label{smo}
\|u\|_{H^{s_0}}\leq c\|Iu\|_{H^{s_{0}+1-s}}\leq cN^{1-s}\|u\|_{H^{s_0}}.
\end{equation}

In the sequel,
we also consider the solutions of the following modified equation
\begin{eqnarray}\label{MNLB}
\left\{
\begin{array}{l}
Iu_{tt}-Iu_{xx}+Iu_{xxxx}-I(|u|^2u)_{xx}=0, \peq  x\in \R,\,t>0,\\
Iu(x,0)=I\phi(x); \quad \partial_tIu(x,0)=I\psi_x(x).
\end{array} \right.
\end{eqnarray}

Applying the operator $(-\Delta)^{-\frac{1}{2}}$ to the equation \eqref{MNLB}, multiplying the result by
$(-\Delta)^{-\frac{1}{2}}\partial_tIu$ and integrating by parts with respect to $x$, we obtain
\begin{equation*}
\frac{1}{2}\frac{d}{dt}\left(\|Iu(t)\|^2_{H^1}+ \|(-\Delta)^{-\frac{1}{2}}\partial_tIu(t)\|^2_{L^2}\right)
+\langle I\left(|u|^2u\right), \partial_tIu\rangle=0.
\end{equation*}

Moreover,
\begin{equation*}
\frac{d}{dt}\|Iu(t)\|^4_{L^4}=4\int_{\R}|Iu|^2Iu\partial_tIu.
\end{equation*}

Therefore
\begin{equation}\label{ACL}
\frac{d}{dt}E(Iu)(t)=\langle |Iu|^2Iu-I\left(|u|^2u\right), \partial_tIu\rangle.
\end{equation}

Most of our arguments here consist in  showing that the quantity $E(Iu)(t)$ is \textit{almost} conserved
in time.

\section{Notations and preliminary results}\label{S3}

Given any positive numbers $a$ and $b$, the notation $a \lesssim b$ means that there exists a positive
constant $\theta$ such that $a \leq \theta b$. Also, we denote $a \sim b$ when, $a \lesssim b$ and
$b \lesssim a$. We use $a+$ and $a-$ to denote $a+\varepsilon$ and $a-\varepsilon$, respectively, for
arbitrarily small exponents $\varepsilon>0$.

To obtain refined estimates in the $X_{s,b}$ spaces, we also use the weighted Sobolev norms given in the
next two definitions.
\begin{defi}\label{BOS+}
For $s,b \in \R$, $X^{+}_{s,b}$ denotes the completion of the Schwartz class $\mathcal{S}(\R^2)$ with
respect to the norm
\begin{equation*}
\|F\|_{X^{+}_{s,b}}=\|\langle\tau+|\xi|^2\rangle^b\langle\xi\rangle^s \widetilde{F}\|_{L^{2}_{\xi,\tau}}.
\end{equation*}
\end{defi}

\begin{defi}\label{BOS-}
For $s,b \in \R$, $X^{-}_{s,b}$ denotes the completion of the Schwartz class $\mathcal{S}(\R^2)$ with respect
to the norm
\begin{equation*}
\|F\|_{X^{-}_{s,b}}=\|\langle\tau-|\xi|^2\rangle^b\langle\xi\rangle^s \widetilde{F}\|_{L^{2}_{\xi,\tau}}.
\end{equation*}
\end{defi}

The following numerical lemma is essential in our arguments.
\begin{lemm}\label{l3.3}
There exists $c>0$ such that
\begin{equation}\label{LN}
\dfrac{1}{c}\leq\sup_{x,y\geq 0}\dfrac{1+|x-y|}{1+|x-\sqrt{y^2+y}|}\leq c.
\end{equation}
\end{lemm}

\noindent\textbf{Proof. } Since $y\leq\sqrt{y^2+y}\leq y+1/2$ for all $y\geq 0$ a simple computation shows the
desired inequalities.\\
\fim

In view of the previous lemma we have an equivalent way to compute the $X_{s,b}$-norm, that is,
\begin{equation*}
\|u\|_{X_{s,b}}\sim\|\langle|\tau|-\xi^2\rangle^b\langle\xi\rangle^{s} \widetilde{u}(\xi,\tau)\|_{L^2_{\xi,\tau}}.
\end{equation*}

Moreover, define $u^{+}$ and $u^{-}$ such that $\widetilde{u^{+}}=\widetilde{u}\chi_{\{\tau\geq 0\}}$
and $\widetilde{u^{-}}=\widetilde{u}\chi_{\{\tau< 0\}}$. Therefore, $\widetilde{u}=\widetilde{u^{+}}
+\widetilde{u^{-}}$ and
\begin{equation*}
\begin{split}
\|\langle|\tau|-\xi^2\rangle^b&\langle\xi\rangle^{s} \widetilde{u}(\xi,\tau)\|^2_{L^2_{\xi,\tau}}=\\
&\|\langle\tau+\xi^2\rangle^b\langle\xi\rangle^{s} \widetilde{u^{+}}(\xi,\tau)\|^2_{L^2_{\xi,\tau}}
+ \|\langle\tau-\xi^2\rangle^b\langle\xi\rangle^{s} \widetilde{u^{-}}(\xi,\tau)\|^2_{L^2_{\xi,\tau}}.
\end{split}
\end{equation*}

In other words, given $u \in X_{s,b}$, there exist $u^{+} \in X^{+}_{s,b}$ and $u^{-} \in X^{-}_{s,b}$
such that
\begin{equation}\label{MXSB}
\|u\|_{X_{s,b}}\sim\|u^{+}\|_{X^{+}_{s,b}}+\|u^{-}\|_{X^{-}_{s,b}}.
\end{equation}

From now on, we refer $u^{+}$ and $u^{-}$ as the decomposition of $u \in X_{s,b}$.

\begin{rema}
Since $X^{+}_{s,b}$ and $X^{-}_{s,b}$ are continuously embedding in $C(\R:H^s(\R))$ for $b>1/2$, we have the same property
for $X_{s,b}$.
\end{rema}

In the context of $X^{\pm}_{s,b}$ spaces, the following two lemmas are well known (see, for example, Ginibre, Tsutsumi and Velo \cite{GTV} and Ozawa and Tsutsumi \cite{OT}).
\begin{lemm}\label{L1}
Let $\psi \in X^{\pm}_{0,\frac{1}{2}+}$, then
\begin{enumerate}
\item [(i)] $\|\psi\|_{L^6_{x,t}}\lesssim \|\psi\|_{X^{\pm}_{0,\frac{1}{2}+}}$,
\item [(ii)]$\|\psi\|_{L^4_{x,t}}\lesssim \|\psi\|_{X^{\pm}_{0,\frac{1}{2}+}}$.
\end{enumerate}
\end{lemm}

\begin{lemm}\label{L2}
Let $\psi_1, \psi_2 \in X^{\pm}_{0,\frac{1}{2}+}$ be supported on spatial frequencies $|\xi|\sim N_1,N_2$,
respectively. Then, if $N_1\lesssim N_2$, one has
\begin{enumerate}
\item [(i)]$\|\psi_1\psi_2\|_{L^2_{x,t}}\lesssim \frac{1}{N_2^{\frac{1}{2}}}\|\psi_1\|_{X^{+}_{0,\frac{1}{2}+}} \|\psi_2\|_{X^{+}_{0,\frac{1}{2}+}}$,
\item [(ii)]$\|\psi_1\psi_2\|_{L^2_{x,t}}\lesssim \frac{1}{N_2^{\frac{1}{2}}}\|\psi_1\|_{X^{-}_{0,\frac{1}{2}+}} \|\psi_2\|_{X^{-}_{0,\frac{1}{2}+}}$,
\item [(iii)] $\|\psi_1\psi_2\|_{L^2_{x,t}}\lesssim \frac{1}{N_2^{\frac{1}{2}}}\|\psi_1\|_{X^{+}_{0,\frac{1}{2}+}} \|\psi_2\|_{X^{-}_{0,\frac{1}{2}+}}.$
\end{enumerate}
\end{lemm}

In view of \eqref{MXSB}, we can also obtain refined Strichartz estimates for the $X_{s,b}$ spaces.

\begin{lemm}\label{L1.2}
Let $\psi \in X^{\delta}_{0,\frac{1}{2}+}$, then
\begin{enumerate}
\item [(i)] $\|\psi\|_{L^6(\R\times[0,\delta])}\lesssim \|\psi\|_{X^{\delta}_{0,\frac{1}{2}+}}$,
\item [(ii)]$\|\psi\|_{L^4(\R\times[0,\delta])}\lesssim \|\psi\|_{X^{\delta}_{0,\frac{1}{2}+}}$.
\end{enumerate}
\end{lemm}

\begin{lemm}\label{L2.2}
Let $\psi_1, \psi_2 \in X^{\delta}_{0,\frac{1}{2}+}$ be supported on spatial frequencies
$|\xi|\sim N_1,N_2$, respectively. Then, if $N_1\lesssim N_2$, one has
\begin{equation}\label{BOUR}
\|\psi_1\psi_2\|_{L^2(\R\times[0,\delta])}\lesssim \frac{1}{N_2^{\frac{1}{2}}}\|\psi_1\|_{X^{\delta}_{0,\frac{1}{2}+}} \|\psi_2\|_{X^{\delta}_{0,\frac{1}{2}+}}.
\end{equation}
\end{lemm}

\noindent\textbf{Proof. }Let $\widetilde{\psi_i}\in X_{0,\frac{1}{2}+}$ be any extension of $\psi_i$ and
($\widetilde{\psi}^{+}_i$, $\widetilde{\psi}^{-}_i$) be its decomposition, for $i=1,2$. Therefore,
we can write
\begin{equation}\label{DECOM}
\widetilde{\psi}_1\widetilde{\psi}_2=\widetilde{\psi}^{+}_1 \widetilde{\psi}^{+}_2 +\widetilde{\psi}^{+}_1\widetilde{\psi}^{-}_2+ \widetilde{\psi}^{-}_1\widetilde{\psi}^{+}_2 + \widetilde{\psi}^{-}_1\widetilde{\psi}^{-}_2.
\end{equation}

Note that this decomposition has an important property: if $\textrm{supp}_x(\widetilde{\psi}_i)\sim N_i$,
then $\textrm{supp}_x(\widetilde{\psi}^{\pm}_i)\sim N_i$. From this fact, we can see that each term on the
right hand side of \eqref{DECOM} satisfies the hypotheses of Lemma \ref{L2}. Therefore, applying the
triangular inequality and the relation \eqref{MXSB} we obtain \eqref{BOUR}.\\
\fim

\section{Local theory}

We first consider the linear equation
\begin{equation}\label{LB}
u_{tt}-u_{xx}+u_{xxxx}=0
\end{equation}
whose solution with initial data $u(0)=\phi$ and $\partial_tu(0)=\psi_x$, is given by
\begin{equation}\label{GUB}
u(t)=V_c(t)\phi+V_s(t)\psi_x,
\end{equation}
where
\begin{eqnarray*}
V_c(t)\phi&=&\left( \frac{e^{it\sqrt{{\xi}^2+{\xi}^4}}+ e^{-it\sqrt{{\xi}^2+{\xi}^4}}}{2}\hat{\phi}(\xi)\right)^{\vee},\\
V_s(t){\psi_x}&=&\left( \frac{e^{it\sqrt{{\xi}^2+{\xi}^4}}- e^{-it\sqrt{{\xi}^2+{\xi}^4}}}{2i\sqrt{{\xi}^2+{\xi}^4}}\hat{\psi_x}(\xi)\right)^{\vee}.
\end{eqnarray*}

By Duhamel's Principle the solution of (NLB) is equivalent to
\begin{equation}\label{INT}
u(t)= V_c(t)\phi+V_s(t)\psi_x+\int_{0}^{t}V_s(t-t')(|u|^2u)_{xx}(t')dt'.
\end{equation}

Let $\theta$ be a cutoff function satisfying $\theta \in C^{\infty}_{0}(\R)$, $0\leq \theta \leq 1$,
$\theta \equiv 1$ in $[-1,1]$, supp$(\theta) \subseteq [-2,2]$ and for $0<\delta<1$ define
$\theta_{\delta}(t)=\theta(t/\delta)$. In fact, to work in the $X_{s,b}$ spaces we consider another
version of (\ref{INT}), that is,
\begin{equation}\label{INT2}
u(t)= \theta(t)\left(V_c(t)\phi+V_s(t)\psi_x\right)+\theta_{\delta}(t)\int_{0}^{t} V_s(t-t')(|u|^2u)_{xx}(t')dt'.
\end{equation}

Note that the integral equation (\ref{INT2}) is defined for all $(x,t)\in \R^2$. Moreover if $u$ is a
solution of \eqref{INT2} then $\tilde{u}=u|_{[0,\delta]}$ will be a solution of (\ref{INT}) in $[0,\delta]$.

In the next two lemmas, we estimate the linear and Duhamel's part of the integral equation (\ref{INT2}).
The proofs can be found in Farah \cite{LG4} Lemmas 2.1-2.2.

\begin{lemm}\label{l21}
Let $u(t)$ the solution of the linear equation
\begin{eqnarray*}
\left\{
\begin{array}{l}
u_{tt}-u_{xx}+u_{xxxx}=0,\\
u(0,x)=\phi(x); \peq \partial_tu(0,x)=(\psi(x))_x
\end{array} \right.
\end{eqnarray*}
with $\phi \in H^s(\R)$ and $\psi \in H^{s-1}(\R)$. Then there exists $c>0$ depending only on $\theta,s,b$
such that
\begin{equation}\label{LP}
\|\theta u\|_{X_{s,b}}\leq c\left(\|\phi\|_{H^s}+\|\psi\|_{H^{s-1}}\right).
\end{equation}
\end{lemm}

%
%
%

\begin{lemm}\label{L22}
Let $-\frac{1}{2}<b'\leq 0\leq b \leq b'+1$ and $0<T \leq 1$, then
\begin{enumerate}
\item [$(i)$] $\left\|\theta_T(t)\int_{0}^{t}g(t')dt'\right\|_{H^b_t} \leq T^{1-(b-b')}\|g\|_{H^{b'}_{t}}$;

\item [$(ii)$] $\left\|\theta_T(t)\int_{0}^{t}V_s(t-t')f(u)(t')dt'\right\|_{X_{s,b}} \leq T^{1-(b-b')} \left\|
f(u)
\right\|_{X_{s,b'}}$.
\end{enumerate}
\end{lemm}

Now, we prove some estimates for the cubic term in (\ref{NLB}).
\begin{lemm}\label{l23}
For $s\geq 0$ and $a\leq 0$ we have that
\begin{itemize}
\item [(i)] $\||u|^2u\|_{X_{s,a}}\lesssim \|u\|^3_{X_{s,\frac{1}{2}+}}$,
\item [(ii)]$\|I(|u|^2u)\|_{X_{1,0}}\lesssim \|Iu\|^3_{X_{1,\frac{1}{2}+}}$.
\end{itemize}
\end{lemm}

\noindent\textbf{Proof. }
\begin{itemize}
\item [(i)] Define $\widetilde{w}(\xi,\tau)=\langle \xi\rangle^s\widetilde{u}(\xi,\tau)$. Since $a\leq 0$
and $\langle \xi\rangle^s\leq \langle \xi-\eta\rangle^s\langle \eta-\sigma\rangle^s\langle \sigma\rangle^s$
for $s\geq 0$, we obtain
\begin{equation}
\begin{split}
\left\||u|^2u\right\|_{X_{s,a}}=&\left\|\langle|\tau|-\gamma(\xi) \rangle^a\langle \xi\rangle^s\widetilde{|u|^2u}(\xi,\tau)\right\|_{L^2}\\
&\leq\sup_{\xi,\tau} \langle|\tau|-\gamma(\xi) \rangle^a\left\||w|^2w\right\|_{L^2}\\
&\leq\left\|w\right\|^3_{L^6}.
\end{split}
\end{equation}
Thus the desired inequality follows by Lemma \ref{L1.2}-($i$).

\item [(ii)] This estimate follows by the previous item together with the interpolation result (Lemma 12.1) of \cite{CKSTT2}.
\end{itemize}
\fim

By standard arguments, the inequalities stated in Lemmas \ref{l21}-\ref{l23} immediately yield the
local well-posedness result stated in Theorem \ref{t3.1}.

Next we prove a local existence result for the equation satisfied by $Iu$ (\ref{MNLB}). Since we do not have
scaling invariance we also need to estimate the solution existence time.

\begin{theo}\label{t3.2}
Assume $s < 1$. Let $(\phi, \psi) \in H^s(\R) \times H^{s-1}(\R)$ be
given. Then there exists a positive number $\delta$ such
that the IVP (\ref{MNLB}) has a unique local solution $Iu\in C([0,\delta]:H^1(\R))$ such that
\begin{equation}\label{CONT}
\|Iu\|_{X^{\delta}_{1,\frac{1}{2}+}}\leq c\left(\|I\phi\|_{H^1}+\|I\psi\|_{L^2}\right).
\end{equation}

Moreover, the existence time can be estimated by
\begin{equation*}
\delta^{\frac{1}{2}-} \sim \dfrac{1}{\left(\|I\phi\|_{H^1}+\|I\psi\|_{L^{2}}\right)^2}.
\end{equation*}
\end{theo}

\noindent\textbf{Proof. } For $({\phi},{\psi})\in H^s(\R)\times H^{s-1}(\R)$, with $s<1$, and $0<\delta\leq 1$
we define the integral equation
\begin{equation}\label{IE}
\Gamma_{\delta}(Iu)(t)= \theta(t)\left(V_c(t)I\phi+V_s(t)I\psi_x\right)+\theta_{\delta}(t)\int_{0}^{t}
V_s(t-t')(I(|u|^2u))_{xx}(t')dt'.
\end{equation}

Our goal is to use the Picard fixed point theorem to find a solution
\begin{equation*}
\Gamma_{\delta}(Iu)=Iu.
\end{equation*}

From \eqref{LP}, Lemma \ref{L22}-($ii$) with $b=\frac{1}{2}+$ and $b'=0$ and Lemma \ref{l23}-($ii$)
we obtain
\begin{eqnarray}\label{C1}
\|\Gamma_{\delta}(Iu)\|_{X_{1,\frac{1}{2}+}}\!\!\!\!&\leq&\!\!\!\! c\left(\|I\phi\|_{H^1}+\|I\psi\|_{L^{2}}\right)+ c\,{\delta}^{\frac{1}{2}-} \left\|Iu\right\|_{X_{1,\frac{1}{2}+}}^3\\
\|\Gamma_{\delta}(Iu)-\Gamma_{\delta}(Iv)\|_{X_{1,\frac{1}{2}+}}\!\!\!\! &\leq&\!\!\!\! c\,{\delta}^{\frac{1}{2}-} \left(\left\|Iu\right\|^2_{X_{1,\frac{1}{2}+}}+ \left\|Iv\right\|^2_{X_{1,\frac{1}{2}+}}\right) \left\|Iu-Iv\right\|_{X_{1,\frac{1}{2}+}}.\nonumber
\end{eqnarray}

We define
\begin{equation*}
X_{1,\frac{1}{2}+}(d)=\left\{u\in X_{1,\frac{1}{2}+}:\|Iu\|_{X_{1,\frac{1}{2}+}}\leq d\right\},
\end{equation*}
where $d=2c\left(\|I\phi\|_{H^1}+\|I\psi\|_{L^{2}}\right)$.\\

Then choosing
\begin{equation*}
0<{\delta}^{\frac{1}{2}-}<\min\left\{\dfrac{1}{4cd^2},1\right\}
\end{equation*}
we have that $\Gamma_{\delta}:X_{1,\frac{1}{2}+}(d)\rightarrow X_{1,\frac{1}{2}+}(d)$ is a contraction and therefore there exists a unique solution ${Iu}\in X_{1,\frac{1}{2}+}(d)$ of (\ref{IE}).

Moreover, we have that $I\tilde{u}=Iu|_{[0,\delta]} \in C([0,\delta]:H^1(\R))\cap X_{1,\frac{1}{2}+}^{\delta}$ is a solution of (\ref{INT}) in $[0,\delta]$.

Finally, estimate (\ref{CONT}) follows by the choice of $\delta$ and inequality (\ref{C1}).\\
\fim

\begin{rema}
For $(\phi, \psi) \in H^s(\R) \times H^{s-1}(\R)$, let $Iu$ be the solution given by
Theorem \ref{t3.2} satisfying the following integral equation
\begin{equation}\label{MINT1}
Iu(t)= V_c(t)I\phi+V_s(t)I\psi_x+\int_{0}^{t} V_s(t-t')(I(|u|^2u))_{xx}(t')dt'
\end{equation}
in the interval $[0,\delta]$.

Thus, a simple calculation implies that
\begin{equation}\label{MINT2}
(-\Delta)^{-\frac{1}{2}}\partial_tIu(t)= G(t)I\phi+H(t)I\psi_x+\int_{0}^{t} H(t-t')(I(|u|^2u))_{xx}(t')dt',
\end{equation}
where $G(\cdot)$ and $H(\cdot)$ denote the multipliers with symbol
$-|\xi|^{-\frac{1}{2}}\gamma(\xi)\sin(\cdot\,\gamma(\xi))$ and
$|\xi|^{-\frac{1}{2}}\cos(\cdot\,\gamma(\xi))$, respectively.

Now, define the function
\begin{equation}\label{MINT3}
v(t)= \theta(t)\left(G(t)I\phi+H(t)I\psi_x\right)+ \theta_{\delta}(t)\int_{0}^{t} H(t-t')(I(|u|^2u))_{xx}(t')dt'.
\end{equation}

Therefore, the definition of $X^{\delta}_{s,b}$ and the same arguments used to prove inequality
\eqref{CONT} implies that
\begin{equation}\label{CONT2}
\|(-\Delta)^{-\frac{1}{2}}\partial_tIu\|_{X^{\delta}_{0,\frac{1}{2}+}}\leq \|v\|_{X_{0,\frac{1}{2}+}}
\leq c\left(\|I\phi\|_{H^1}+\|I\psi\|_{L^2}\right).
\end{equation}
\end{rema}

\section{Almost conservation law}

In this section we will establish estimates to control the growth rate of $E(Iu)(t)$.

\begin{prop}\label{p4.1}
Let $s>1/2$, $N\gg 1$ and $Iu$ be a solution of \eqref{MNLB} on $[0, \delta]$ in the sense of Theorem \ref{t3.2}. Then the following estimate holds
\begin{equation}\label{CC}
\left|E(Iu)(\delta)-E(Iu)(0)\right|\leq N^{-2+}\left\|Iu\right\|_{X^{\delta}_{1,\frac{1}{2}+}}^3 \left\|(-\Delta)^{-\frac{1}{2}}\partial_tIu\right\|_{X^{\delta}_{0,\frac{1}{2}+}}.
\end{equation}
\end{prop}

\begin{rema}
The exponent $-2+$ on the right hand side of (\ref{CC}) is directly tied to the restriction $s > 2/3$ in our main theorem. If one could replace the increment $N^{-2+}$ by $N^{-\alpha+}$ for some $\alpha>0$ the argument we give in Section $6$ would imply global well-posedness of (\ref{NLB}) for all $s > 1-\alpha/6$. 
\end{rema}

\noindent\textbf{Proof. } Following the arguments used in \cite{CKSTT4} we only need to bound
\begin{equation*}
E(Iu)(\delta)-E(Iu)(0)=
\end{equation*}
\begin{equation*}
\int_0^{\delta}\!\!\!\int_{\sum_{i=1}^4\xi_i=0} \left( 1-\frac{m(\xi_2+\xi_3+\xi_4)}{m(\xi_2)m(\xi_3)m(\xi_4)}\right) \widehat{{\partial_tIu(\xi_1)}}
 \widehat{Iu(\xi_2)}
 \widehat{{Iu(\xi_3)}}
 \widehat{u(\xi_4)}.
\end{equation*}

Therefore, our aim is to obtain the following inequality
\begin{equation*}
\mathbf{Term}\leq N^{-2+}\left\|I\phi_1\right\|_{X^{\delta}_{0,\frac{1}{2}+}} \prod_{i=2}^4\left\|I\phi_i\right\|_{X^{\delta}_{1,\frac{1}{2}+}},
\end{equation*}
where
\begin{equation*}
\mathbf{Term}\equiv \left|\int_0^{\delta}\!\!\!\int_{\ast} \left( 1-\frac{m(\xi_2+\xi_3+\xi_4)}{m(\xi_2)m(\xi_3)m(\xi_4)}\right)|\xi_1| \widehat{{I\phi_1(\xi_1)}}
 \widehat{I\phi_2(\xi_2)}
 \widehat{{I\phi_3(\xi_3)}}
 \widehat{I\phi_4(\xi_4)}\right|
\end{equation*}
and $\ast$ denotes integration over $\sum_{i=1}^4\xi_i=0$.

We estimate $\mathbf{Term}$ as follows. Without loss of generality, we assume the Fourier transforms of
all these functions to be nonnegative. First, we bound the symbol in the parentheses pointwise
in absolute value, according to the relative sizes of the frequencies involved. After that, the
remaining integrals are estimated using Plancherel formula, H\"older's inequality and Lemmas
\ref{L1.2}-\ref{L2.2}. To sum over the dyadic pieces at the end we need to have extra factors
$N_j^{0-}$, $j=1,2,3,4,$, everywhere.

We decompose the frequencies $\xi_j$, $j=1,2,3,4,$ into dyadic blocks $N_j$. By the symmetry
of the multiplier
\begin{equation}\label{MULT}
1-\frac{m(\xi_2+\xi_3+\xi_4)}{m(\xi_2)m(\xi_3)m(\xi_4)}
\end{equation}
in $\xi_2$, $\xi_3$, $\xi_4$, we may assume that
\begin{equation*}
N_2\geq N_3\geq N_4.
\end{equation*}

Moreover, we can assume
$N_2 \gtrsim N$, because otherwise the symbol is zero. The condition $\sum_{i=1}^{4}\xi_i=0$
implies $N_1\lesssim N_2$. We split the different frequency interaction into several cases,
according to the size of the parameter $N$ in comparison to the $N_i$'s.\\

\quebra \textbf{Case $A$: }$N_2\gtrsim N\gg N_3\geq N_4$.\\

The condition $\sum_{i=1}^{4}\xi_i=0$ implies $N_1\sim N_2$. By mean value theorem,
\begin{equation*}
\left|\frac{m(\xi_2)-m(\xi_2+\xi_3+\xi_4)}{m(\xi_2)}\right|\lesssim \frac{\left|\nabla m(\xi_2)(\xi_3+\xi_4)\right|}{m(\xi_2)}\lesssim\frac{N_3}{N_2}.
\end{equation*}

Therefore, Lemma \ref{L2.2} implies that
\begin{eqnarray*}
{\mathbf{Term}} \!\!\!&\lesssim& \!\!\! \frac{N_1N_3}{N_2}\left\|I\phi_1I\phi_3\right\|_{L^2(\R\times[0,\delta])} \left\|I\phi_2I\phi_4\right\|_{L^2(\R\times[0,\delta])}\\
&\lesssim&  \!\!\!\frac{N_1N_3}{N_2N_1^{\frac{1}{2}}N_2^{\frac{1}{2}}N_2\langle N_3 \rangle\langle N_4 \rangle}\|I\phi_1\|_{X^{\delta}_{0,\frac{1}{2}+}} \prod_{i=2}^{4} \|I\phi_i\|_{X^{\delta}_{1,\frac{1}{2}+}}\\
&\lesssim& \!\!\! N^{-2+}N_{max}^{0-}\|I\phi_1\|_{X^{\delta}_{0,\frac{1}{2}+}} \prod_{i=2}^{4} \|I\phi_i\|_{X^{\delta}_{1,\frac{1}{2}+}}.
\end{eqnarray*}
\textbf{Case $B$: }$N_2\gg N_3\gtrsim N$ and $N_3\geq N_4$.\\

In this case we also have $N_1\sim N_2$. We bound the multiplier (\ref{MULT}) by
\begin{equation}\label{MULT2}
\left|1-\frac{m(\xi_2+\xi_3+\xi_4)}{m(\xi_2)m(\xi_3)m(\xi_4)}\right|\lesssim \frac{m(\xi_1)}{m(\xi_2)m(\xi_3)m(\xi_4)}.
\end{equation}

Therefore, since $m(N_1)\sim m(N_2)$, applying Lemma \ref{L2.2} we have
\begin{eqnarray*}
{\mathbf{Term}} \!\!\!&\lesssim&  \!\!\!\frac{N_1}{m(N_3)m(N_4)}\left\|I\phi_1I\phi_3\right\|_{L^2(\R\times[0,\delta])} \left\|I\phi_2I\phi_4\right\|_{L^2(\R\times[0,\delta])}\\
&\lesssim& \!\!\! \frac{N_1} {m(N_3)m(N_4)N_1^{\frac{1}{2}}N_2^{\frac{1}{2}}N_2N_3\langle N_4 \rangle}\|I\phi_1\|_{X^{\delta}_{0,\frac{1}{2}+}} \prod_{i=2}^{4} \|I\phi_i\|_{X^{\delta}_{1,\frac{1}{2}+}}\\
&\lesssim& \!\!\! \frac{1} {m(N_3)N_3m(N_4)\langle N_4 \rangle N_2}\|I\phi_1\|_{X^{\delta}_{0,\frac{1}{2}+}} \prod_{i=2}^{4} \|I\phi_i\|_{X^{\delta}_{1,\frac{1}{2}+}}\\
&\lesssim& \!\!\! N^{-2+}N_{max}^{0-}\|I\phi_1\|_{X^{\delta}_{0,\frac{1}{2}+}} \prod_{i=2}^{4} \|I\phi_i\|_{X^{\delta}_{1,\frac{1}{2}+}},
\end{eqnarray*}
where in the last inequality we use the fact that for any $p \geq \frac{1}{2}$, the function
$m(x)x^p$ is increasing and $m(x)\langle x\rangle^p$ is bounded below, which implies
$m(N_3)N_3\gtrsim m(N)N=N$ and $m(N_4)\langle N_4\rangle \gtrsim 1$.\\

\quebra \textbf{Case $C$: }$N_2\sim N_3\gtrsim N$ and $N_3\geq N_4$.\\

The condition $\sum_{i=1}^{4}\xi_i=0$ implies $N_1\lesssim N_2$. We again bound the multiplier (\ref{MULT}) pointwise by (\ref{MULT2}). To obtain the decay $N^{-2+}$ we split this case into four subcases.\\

\textbf{Case $C.1$: } $N_4\gtrsim N$ and $N_4\ll N_3$.\\

From (\ref{MULT2}) and Lemmas \ref{L1.2}-\ref{L2.2}, we have that
\begin{equation*}
\begin{split}
{\mathbf{Term}} &\lesssim \frac{N_1m(N_1)}{m(N_2)m(N_3)m(N_4)} \prod_{i=\{1,3\}}\left\|I\phi_i\right\|_{L^4(\R\times[0,\delta])} \left\|I\phi_2I\phi_4\right\|_{L^2(\R\times[0,\delta])}\\
&\lesssim \frac{N_1m(N_1)} {m(N_2)m(N_3)m(N_4)N_2^{\frac{1}{2}}N_2N_3 N_4 }\|I\phi_1\|_{X^{\delta}_{0,\frac{1}{2}+}} \prod_{i=2}^{4} \|I\phi_i\|_{X^{\delta}_{1,\frac{1}{2}+}}
\end{split}
\end{equation*}
\begin{equation*}
\begin{split}
&\lesssim \! \frac{N_{max}^{0-}} {m(N_2)N_2^{\frac{3}{4}-}m(N_3)N_3^{\frac{3}{4}}m(N_4) N_4}\|I\phi_1\|_{X^{\delta}_{0,\frac{1}{2}+}} \prod_{i=2}^{4} \|I\phi_i\|_{X^{\delta}_{1,\frac{1}{2}+}}\\
&\lesssim \! N^{-2+}N_{max}^{0-}\|I\phi_1\|_{X^{\delta}_{0,\frac{1}{2}+}} \prod_{i=2}^{4} \|I\phi_i\|_{X^{\delta}_{1,\frac{1}{2}+}}.
\end{split}
\end{equation*}

\textbf{Case $C.2$: } $N_4\gtrsim N$ and $N_4\sim N_3$.\\

Applying the same arguments as above
\begin{eqnarray*}
{\mathbf{Term}} \!\!\!&\lesssim& \!\!\! \frac{N_1m(N_1)}{m(N_2)m(N_3)m(N_4)}  \prod_{i=1}^4\left\|I\phi_i\right\|_{L^4(\R\times[0,\delta])} \\
&\lesssim&  \!\!\!\frac{N_1m(N_1)} {m(N_2)m(N_3)m(N_4) N_2N_3 N_4 }\|I\phi_1\|_{X^{\delta}_{0,\frac{1}{2}+}} \prod_{i=2}^{4} \|I\phi_i\|_{X^{\delta}_{1,\frac{1}{2}+}}\\
&\lesssim&  \!\!\!\frac{N_{max}^{0-}} {m(N_2)N_2^{\frac{2}{3}-}m(N_3)N_3^{\frac{2}{3}}m(N_4) N_4^{\frac{2}{3}}}\|I\phi_1\|_{X^{\delta}_{0,\frac{1}{2}+}} \prod_{i=2}^{4} \|I\phi_i\|_{X^{\delta}_{1,\frac{1}{2}+}}\\
&\lesssim& \!\!\! N^{-2+}N_{max}^{0-}\|I\phi_1\|_{X^{\delta}_{0,\frac{1}{2}+}} \prod_{i=2}^{4} \|I\phi_i\|_{X^{\delta}_{1,\frac{1}{2}+}}.
\end{eqnarray*}

\textbf{Case $C.3$:} $N_4\ll N$ and $N_1\ll N_2$.\\

Again using the bound (\ref{MULT2}) and Lemma \ref{L2.2}, we have
\begin{eqnarray*}
{\mathbf{Term}} \!\!\!&\lesssim& \!\!\! \frac{N_1m(N_1)}{m(N_2)m(N_3)m(N_4)} \left\|{I\phi_1}I\phi_2\right\|_{L^2(\R\times[0,\delta])} \left\|{I\phi_3}I\phi_4\right\|_{L^2(\R\times[0,\delta])}\\
&\lesssim& \!\!\! \frac{N_1m(N_1)} {m(N_2)m(N_3)m(N_4)N_2^{\frac{1}{2}}N_3^{\frac{1}{2}} N_2N_3 \langle N_4\rangle }\|I\phi_1\|_{X^{\delta}_{0,\frac{1}{2}+}} \prod_{i=2}^{4} \|I\phi_i\|_{X^{\delta}_{1,\frac{1}{2}+}}\\
&\lesssim&  \!\!\!\frac{N_{max}^{0-}} {m(N_2)N_2^{1-}m(N_3)N_3m(N_4) \langle N_4\rangle}\|I\phi_1\|_{X^{\delta}_{0,\frac{1}{2}+}} \prod_{i=2}^{4} \|I\phi_i\|_{X^{\delta}_{1,\frac{1}{2}+}}\\
&\lesssim&  \!\!\!N^{-2+}N_{max}^{0-}\|I\phi_1\|_{X^{\delta}_{0,\frac{1}{2}+}} \prod_{i=2}^{4} \|I\phi_i\|_{X^{\delta}_{1,\frac{1}{2}+}}.
\end{eqnarray*}

\textbf{Case $C.4$: } $N_4\ll N$ and $N_1\sim N_2\sim N_3 \gtrsim N$.\\

In this case, we use an argument similar to the one used in Pecher \cite{P1} Proposition 5.1. Because
of $\sum_{i=1}^{4}\xi_i=0$, two of the large frequencies have different sign,
say, $\xi_1$ and $\xi_2$. Thus,
\begin{equation*}
|\xi_1|\leq |\xi_1-\xi_2|\leq 2|\xi_1|
\end{equation*}
and
\begin{equation*}
|\xi_1+\xi_2|= |\xi_3+\xi_4|\sim|\xi_1|.
\end{equation*}

Therefore, using the bound \eqref{MULT2} and Lemma \ref{L2.2}, we have
\begin{eqnarray*}
{\mathbf{Term}} \!\!\!&\lesssim& \!\!\!\frac{N_1^{\frac{1}{2}}m(N_1)}{m(N_2)m(N_3)m(N_4)} \left\|(D_x^{\frac{1}{2}}I\phi_1)I\phi_2\right\|_{L^2_{x,t}} \left\|I\phi_3I\phi_4\right\|_{L^2_{x,t}}\\
&\lesssim& \!\!\! \frac{N_1^{\frac{1}{2}}m(N_1)} {m(N_2)m(N_3)m(N_4)N_3^{\frac{1}{2}} N_3 \langle N_4\rangle }\left\|(D_x^{\frac{1}{2}}I\phi_1)I\phi_2\right\|_{L^2_{x,t}} \prod_{i=3}^{4} \|I\phi_i\|_{X^{\delta}_{1,\frac{1}{2}+}}\\
&\lesssim& \!\!\! \frac{N_{max}^{0-}} {m(N_2)N_2^{1-}m(N_3)N_3m(N_4) \langle N_4\rangle}\|I\phi_1\|_{X^{\delta}_{0,\frac{1}{2}+}} \prod_{i=2}^{4} \|I\phi_i\|_{X^{\delta}_{1,\frac{1}{2}+}}\\
&\lesssim& \!\!\! N^{-2+}N_{max}^{0-}\|I\phi_1\|_{X^{\delta}_{0,\frac{1}{2}+}} \prod_{i=2}^{4} \|I\phi_i\|_{X^{\delta}_{1,\frac{1}{2}+}},
\end{eqnarray*}
where we have estimated $\|(D_x^{\frac{1}{2}}I\phi_1)I\phi_2\|_{L^2_{x,t}}$ via the following lemma.

\begin{lemm}\label{L3}
Let $\psi_1, \psi_2 \in X_{0,\frac{1}{2}+}$ be supported on spatial frequencies $|\xi_i|\sim N_i$, $i=1,2$.
If $|\xi_1|\lesssim \min\left\{|\xi_1-\xi_2|,|\xi_1+\xi_2|\right\}$ for all
$\xi_i\in \textrm{supp}(\widehat{\psi}_i)$, $i=1,2$, then
\begin{equation}\label{GRU}
\|(D^{\frac{1}{2}}_x\psi_1)\psi_2\|_{L^2_{x,t}}\lesssim \|\psi_1\|_{X^{\delta}_{0,\frac{1}{2}+}}
 \|\psi_2\|_{X^{\delta}_{0,\frac{1}{2}+}}.
\end{equation}
\end{lemm}

\noindent\textbf{Proof of Lemma \ref{L3}:} By the same arguments used in the proof of Lemma \ref{L2.2} we just
need to establish  the following inequalities
\begin{enumerate}
\item [(i)] $\|(D^{\frac{1}{2}}_x\psi_1^{\pm})\psi_2^{\pm}\|_{L^2_{x,t}}\lesssim \|\psi_1^{\pm}\|_{X^{\pm}_{0,\frac{1}{2}+}} \|\psi_2^{\pm}\|_{X^{\pm}_{0,\frac{1}{2}+}}$
\item [(ii)]$\|(D^{\frac{1}{2}}_x\psi_1^{\pm})\psi_2^{\mp}\|_{L^2_{x,t}}\lesssim \|\psi_1^{\pm}\|_{X^{\pm}_{0,\frac{1}{2}+}} \|\psi_2^{\mp}\|_{X^{\mp}_{0,\frac{1}{2}+}}$
\end{enumerate}
for all $\psi_1^{\pm}, \psi_2^{\pm} \in X^{\pm}_{0,\frac{1}{2}+}$ supported on spatial frequencies $|\xi_i|\sim N_i$, $i=1,2$, such that $|\xi_1|\lesssim \min\left\{|\xi_1-\xi_2|,|\xi_1+\xi_2|\right\}$.\\

\underline{ Proof of (i) }\\

Applying standard arguments, it is sufficient to prove that
\begin{equation*}
\|(e^{{\mp}it\Delta} D^{\frac{1}{2}}_xu^{\pm})(e^{{\mp}it\Delta} v^{\pm})\|_{L^2_{x,t}}\lesssim \|u^{\pm}\|_{L^2} \|v^{\pm}\|_{L^2}
\end{equation*}
for functions $\widehat{u}^{\pm}(\xi_1)$ and $\widehat{v}^{\pm}(\xi_2)$ with support in $|\xi_i|\sim N_i$, $i=1,2$.

We have
\begin{eqnarray*}
\|(e^{{\mp}it\Delta} D^{\frac{1}{2}}_xu^{\pm})(e^{{\mp}it\Delta} v^{\pm})\|_{L^2_{x,t}}
\end{eqnarray*}
\vspace*{-0.25in}
\begin{eqnarray*}
&=&\int\!\!\!\int\left(\int\!\!\!\int_{\ast} e^{{\mp}it(\xi_1^2+\xi_2^2-\eta_1^2-\eta_2^2)} |\xi_1|\widehat{u}^{\pm}(\xi_1)\widehat{v}^{\pm}(\xi_2) \widehat{u}^{\pm}(\eta_1)\widehat{v}^{\pm}(\eta_2)d\xi_1d\eta_1\right)d\xi dt\\
&=&\int\left(\int\!\!\!\int_{\ast} \delta(P_{\pm}(\eta_1))|\xi_1|\widehat{u}^{\pm}(\xi_1)\widehat{v}^{\pm}(\xi_2) \widehat{u}^{\pm}(\eta_1)\widehat{v}^{\pm}(\eta_2)d\xi_1d\eta_1\right)d\xi
\end{eqnarray*}
where $\ast$ denotes integration over $\xi=\xi_1+\xi_2=\eta_1+\eta_2$ and
\begin{equation*}
P_{\pm}(\eta_1)=\pm(\xi_1^2+\xi_2^2-\eta_1^2-\eta_2^2).
\end{equation*}

Note that $P_{\pm}(\eta_1)$ has roots $\eta_1=\xi_1$ and $\eta_1=\xi-\xi_1$. Now, using the
well-known identity $\delta(g(x))=\sum_n\frac{\delta(x-x_n)}{|g^{'}(x_n)|}$, where the sum is
taken over all simple zeros of $g$, we obtain
\begin{eqnarray*}
\|(e^{{\mp}it\Delta} D^{\frac{1}{2}}_xu^{\pm})(e^{{\mp}it\Delta} v^{\pm})\|_{L^2_{x,t}}
\end{eqnarray*}
\vspace*{-0.25in}
\begin{eqnarray*}
&\lesssim&\int\left(\int \dfrac{|\xi_1|\widehat{u}^{\pm}(\xi_1) \widehat{u}^{\pm}(\xi_1)\widehat{v}^{\pm}(\xi-\xi_1) \widehat{v}^{\pm}(\xi-\xi_1)}{|\xi_1-\xi_2|} d\xi_1\right)d\xi\\
&&+\int\left(\int \dfrac{|\xi_1|\widehat{u}^{\pm}(\xi_1) \widehat{u}^{\pm}(\xi-\xi_1)\widehat{v}^{\pm}(\xi-\xi_1) \widehat{v}^{\pm}(\xi_1)}{|\xi_1-\xi_2|} d\xi_1\right)d\xi\\
&\lesssim&\|u^{\pm}\|_{L^2} \|v^{\pm}\|_{L^2},
\end{eqnarray*}
where in the last inequality we have used the fact that $|\xi_1|\leq |\xi_1-\xi_2|$.\\

\underline { Proof of (ii) }\\

Again it is sufficient to prove
\begin{equation*}
\|(e^{{\mp}it\Delta} D^{\frac{1}{2}}_xu^{\pm})(e^{{\pm}it\Delta} v^{\mp})\|_{L^2_{x,t}}
\lesssim \|u^{\pm}\|_{L^2} \|v^{\pm}\|_{L^2}
\end{equation*}
for functions $\widehat{u}^{\pm}(\xi_1)$ and $\widehat{v}^{\pm}(\xi_2)$ with support in
$|\xi_i|\sim N_i$, $i=1,2$.

We have
\begin{eqnarray*}
\|(e^{{\mp}it\Delta} D^{\frac{1}{2}}_xu^{\pm})(e^{{\pm}it\Delta} v^{\mp})\|_{L^2_{x,t}}
\end{eqnarray*}
\vspace*{-0.25in}
\begin{eqnarray*}
&=&\int\left(\int\!\!\!\int_{\ast} \delta(P_{\pm}(\eta_1))|\xi_1|\widehat{u}^{\pm}(\xi_1)\widehat{v}^{\mp}(\xi_2) \widehat{u}^{\pm}(\eta_1)\widehat{v}^{\mp}(\eta_2)d\xi_1d\eta_1\right)d\xi
\end{eqnarray*}
where $\ast$ denotes integration over $\xi=\xi_1+\xi_2=\eta_1+\eta_2$ and
\begin{equation*}
P_{\pm}(\eta_1)=\mp(\xi_1^2-\xi_2^2-\eta_1^2+\eta_2^2).
\end{equation*}

Since $P_{\pm}(\eta_1)$ has root $\eta_1=\xi_1$ and $P^{'}_{\pm}(\xi_1)=2|\xi_1+\xi_2|$, we obtain
\begin{equation*}
\begin{split}
\|(e^{{\mp}it\Delta} D^{\frac{1}{2}}_x&u^{\pm})(e^{{\pm}it\Delta} v^{\mp})\|_{L^2_{x,t}}\\
&\lesssim\int\left(\int \dfrac{|\xi_1|\,\widehat{u}^{\pm}(\xi_1) \widehat{u}^{\pm}(\xi_1)\widehat{v}^{\mp}(\xi-\xi_1) \widehat{v}^{\mp}(\xi-\xi_1)}{|\xi_1+\xi_2|}\, d\xi_1\right)d\xi\\
&\lesssim\|u^{\pm}\|_{L^2} \|v^{\pm}\|_{L^2},
\end{split}
\end{equation*}
where in the last inequality we use the fact that $|\xi_1|\leq |\xi_1+\xi_2|$.\\

This completes the proof of Proposition \ref{p4.1}.\\
\fim

\section{Global theory}

Once the relation ($\ref{CC}$) is obtained, we can proof our global result stated in Theorem \ref{T1}.\\

\noindent\textbf{Proof of Theorem \ref{T1}. }Let $(\phi, \psi) \in H^s(\R) \times H^{s-1}(\R)$
with $1/2\leq s <1$. From the definition of the multiplier $I$ we have
\begin{eqnarray*}
\|I\phi\|_{H^1}&\leq& c N^{1-s}\|\phi\|_{H^{s}},\\
\|I\psi\|_{L^2}&\leq& c N^{1-s}\|\psi\|_{H^{s-1}}.
\end{eqnarray*}

Therefore, there exists $c_1>0$ such that
\begin{equation*}
\max\left\{\|I\phi\|_{H^1}, \|I\psi\|_{L^2}\right\}\leq c_1 N^{1-s}.
\end{equation*}

We use our local existence theorem on $[0, \delta]$, where $\delta^{\frac{1}{2}-}\sim N^{-2(s-1)}$ and conclude
\begin{eqnarray*}
\max\left\{\|Iu\|_{X^{\delta}_{1,\frac{1}{2}+}}, \|(-\Delta)^{-\frac{1}{2}}\partial_tIu\|_{X^{\delta}_{0,\frac{1}{2}+}} \right\}&\leq& c\left(\|I\phi\|_{H^1}+\|I\psi\|_{L^2}\right)\\
&\leq& c_2 N^{1-s}.
\end{eqnarray*}

From the conservation law (\ref{EC}), we obtain
\begin{equation}\label{CC2}
\|Iu(\delta)\|^2_{H^1}+\|(-\Delta)^{-\frac{1}{2}}\partial_tIu(\delta)\|^2_{L^2}\leq c_3 E(Iu(\delta)).
\end{equation}

On the other hand, since $\|f\|_{L^4}\lesssim \|f\|_{H^s}$ for $s\geq 1/4$, we have
\begin{equation*}
E(Iu(0))\leq c \,N^{2(1-s)}+c\|\phi\|_{L^4}^4\leq c_4N^{2(1-s)}.
\end{equation*}

By the almost conservation law stated in Proposition \ref{p4.1}, we have
\begin{equation*}
\left|E(Iu)(\delta)-E(Iu)(0)\right|\leq c_5N^{-2+}\left\|Iu\right\|_{X^{\delta}_{1,\frac{1}{2}+}}^3 \left\|(-\Delta)^{-\frac{1}{2}}\partial_tIu\right\|_{X^{\delta}_{0,\frac{1}{2}+}}.
\end{equation*}

Given a time $T>0$, the number of iteration steps to reach this time is $T\delta^{-1}$. To reapply
the local existence result with time intervals of equal length we need a uniform bound of the solution
at time $t = \delta$ and $t = 2\delta$ etc. In view of (\ref{CC2}), this uniform bound can be obtained
if we control the growth of $E(Iu)(\cdot)$ on the interval $[0,T]$. Therefore, to carry out $T\delta^{-1}$
iterations on time intervals, before the quantity $E(Iu)(t)$ doubles, the following condition has to be
fulfilled
\begin{eqnarray}\label{GGV2}
N^{-2+}N^{4(1-s)}T\delta^{-1}\ll N^{2(1-s)}.
\end{eqnarray}

Since $\delta^{\frac{1}{2}-}\sim N^{-2(1-s)}$, the condition (\ref{GGV2}) can be obtained for 
\begin{equation}\label{expot}
T \sim N^{(6s-4)-}.
\end{equation}

\begin{rema}
Note that the exponent of $N$ on the right hand side of (\ref{expot}) is positive provided $s>2/3$, hence the definition of $N$ makes sense for arbitrary large $T$.
\end{rema}

Therefore, by our choice of $N$, relation (\ref{smo}) and (\ref{CC2}) imply for $T\gg 1$ that
\begin{equation*}
\sup_{t\in[0,T]}\left\{\|u(t)\|^2_{H^{s}}+ \|(-\Delta)^{-\frac{1}{2}}\partial_tu(t)\|^2_{H^{s-1}}\right\}\lesssim E(Iu(T)) \lesssim T^{\frac{2(1-s)}{6s-4}+}
\end{equation*}
which implies the polynomial bound (\ref{pb}).\\
\fim



E-mail: farah@impa.br and linares@impa.br.
\end{document}